\documentclass[12pt]{article}
\newcommand{\ba}{\noindent $\begin{array}}
\newcommand{\ea}{\end{array}$}
\newcommand{\be}{\begin{equation}}
\newcommand{\ee}{\end{equation}}
\newcommand{\bd}{\begin{displaymath}}
\newcommand{\ed}{\end{displaymath}}
\newcommand{\beq}{\begin{eqnarray*}}
\newcommand{\eeq}{\end{eqnarray*}}
\newcommand{\beqn}{\begin{eqnarray}}
\newcommand{\eeqn}{\end{eqnarray}}
\usepackage{amsfonts,latexsym,lscape}

\small\normalsize

\usepackage{rawfonts}
\newtheorem{theorem}{Theorem}[section]

\newtheorem{definition}{Definition}[section]

\newfont{\Bb}{msbm10 scaled\magstep1}

\begin{document}

\pagestyle{plain} \thispagestyle{empty}
\begin{center}
               { \textbf{Special Report for \\
         3rd INTERNATIONAL  CONFERENCE ON THE ABS ALGORITHMS\\
                    May 13-14/01, Beijing}}
\end{center}
\vspace{1.5cm}
\begin{center}
{

      \large \textbf{AN ABS ALGORITHM FOR A CLASS OF SYSTEMS OF STOCHASTIC LINEAR EQUATIONS}
\vspace{2cm}

                  Hai-Shan Han, Zun-Quan Xia and Antonino Del Popolo\\
\vspace{1cm}
 \vspace{10cm}
                   Centre for Optimization Research and Applications\\
        Department of Applied Mathematics, Dalian University of Technology
        }
\end{center}

\newpage

\begin{center}
{\large \textbf{AN ABS ALGORITHM FOR A CLASS OF SYSTEMS OF
STOCHASTIC LINEAR EQUATIONS}}
\end{center}

\begin{center}
{\textbf{
                  Hai-Shan Han\footnote{CORA, Department of Applied Mathematics,
                               Dalian University of Technology
                               Dalian 116024, China; Department of Mathematics,Inner Mongolia University of Nationalities,Tongliao 028000 China},
                 Zun-Quan Xia\footnote{CORA, Department of Applied Mathematics,
                               Dalian University of Technology
                               Dalian 116024, China, zqxiazhh@dlut.edu.cn}
                 Antonino Del Popolo\footnote{Dept. of Mathematics, University of Bergamo,
                               Piazza Rosate 2, 24129 Bergamo, Italy}
                }
}
\end{center}

\vspace{2mm} \vspace{2mm}
\begin{center}
\parbox{13.5cm}{\small
\textbf{Abstract.}  This paper is to explore a model of the ABS
Algorithms for dealing with a class of systems of linear
stochastic equations $A\xi=\eta$ satisfying $\eta \sim N_m(v,
I_{m})$. It is shown that the iteration step $\alpha_{i}$ is
$N(V,\pi)$ and approximation solutions is $\xi_{i} \sim
N_n(U,\Sigma)$ for this algorithm model. And some properties of
$(V,\pi)$ and $(U,\Sigma)$ are given.
\\[10pt]
                           }
\end{center}

\section{Introduction}
It is well known, that in stochastic programming, stochastic
linear equations being of the form
\begin{equation}\label{hhx1}
   A(w)\xi(w)=\eta(w)
\end{equation}
or
\begin{equation}\label{hhx2}
  A\xi(w)=\eta(w)
\end{equation}
  with $\xi(w)$\, $\eta(w)$ are random vectors. If
the constraint equation of the mathematical programming, $Ax=b$
contain the random variables, then the stochastic linear equation
has the form (\ref{hhx1}) or (\ref{hhx2}).

In this paper we try explore an application of the classical ABS
algorithms to solve system indicated by (\ref{hhx2}), in other
words, we attempt to establish an model of algorithm, called
 ABS-S defined in Section 3 for solving (\ref{hhx2}). To this end, it is
necessary to recall the ABS algorithms.

Consider the general linear systems, where rank($A$)  is
arbitrary,
\begin{equation}\label{hhx3}
    Ax=b
\end{equation}
or
\begin{equation}\label{hhx4}
    a_i^Tx=b_i,\quad i=1,\cdots,m
\end{equation}
where $x \in I\!\!R^m$,\quad $b \in I\!\!R^n$,\quad $m \leq n$
\quad and \quad$A \in I\!\!R^{m \times n}$, $A=(a_1, a_2,\cdots,
a_m )^T$.

The class of ABS algorithms originally for solving (\ref{hhx3})
or (\ref{hhx4}) was introduced by [AbBs84] and [AbSp89].
 The iterate  scheme
of the basic ABS class of algorithm is defined as follows:

\bigskip
\noindent \textbf{Basic ABS Class of Algorithms}:
see,\,[AbBs84]\,,\,\, [AbBs89]:
\begin{description}
\item[\textbf {  (A) }]  Initialization.\\[1pt]
                     Give  an arbitrary vector $ x_1 \in  I\!\!R^n $,  and an arbitrarily nonsingular matrix $ H_1 \in R^{n,n} $.\\
                     Set $ i = 1 $ and iflag=0.
\vspace{-2mm}
\item[\textbf{  (B)}] Computer two quantities.\\[1pt]
                      Compute
                       $$
                          \begin{array}{l}
                                   s_i = H_i a_i \\
                                   \tau_i=\tau^Te_i =a_i^T x_i-b^T e_i
                          \end{array}
                       $$
\vspace{-2mm}
\item[\textbf{  (C)}] Check the compitability of the system of linear equations.\\[1pt]
                    If $ s_i \not= 0 $ then goto (D). \\
                    If $ s_i = 0 $ and $ \tau_i = 0 $ then set
                      $$
                                    \begin{array}{l}
                                         x_{i+1} = x_i \\
                                         H_{i+1} = H _i
                                  \end{array}
                     $$
                  and goto (F), the $ i $-th equation  is a linear combination of the previous equations. Otherwise  stop, the
                  system has no solution.
\vspace{-2mm}
\item[\textbf{ (D)}] Computer the search vector $p_i\in I\!\!R^n$ by
                   \begin{equation}\label{F2.1}
                        p_i=H_i^Tz_i
                   \end{equation}
                  where $z_{i}$, the  parameter of Broyden, is arbitrary
                  save that
                   \begin{equation}\label{F2.2}
                          z_i^TH_ia_i\neq 0
                    \end{equation}
\vspace{-2mm}
\item[\textbf{ (E)}] Update the approximation of  the solution $ x_i $ by
                   \begin{equation}\label{F2.3}
                        x_{i+1} = x_i - \alpha _i p_i
                    \end{equation}
                where the stepsize $ \alpha _i $ is computed by
                   \begin{equation}\label{F2.4}
                        \alpha _i = \tau_i /a_i^T p_i
                   \end{equation}
                  If $i=m$ stop; $x_{m+1}$ solves the system.
\vspace{-2mm}
\item[\textbf{ (F)}] Update the (Abaffian) matrix $ H_i. $
                   Compute
                   \begin{equation}\label{F2.5}
                        H_{i+1} = H_i - H_ia_iw_i^TH_i /w_i^TH_ia_i
                   \end{equation}
                   where $ w_i \in I\!\!R^n $, the parameter of Abaffy, is arbitrary save for the condition
                    \begin{equation}\label{F2.6}
                              w_i^T H_i a_i = 1 \textrm{ or } \neq 0
                     \end{equation}
\vspace{-2mm}
\item[\textbf{ (G)}] Increment the index $ i $ by one and goto (B).
\end{description}

\bigskip
\noindent We define $n$ by $i$ matrices $A_i,\, W_i$ and $ P_i$ by
\begin{equation}\label{F2.7}
    \begin{array}{l}
                  A_i=(a_1,\,\cdots,\,a_i)^T, \quad
                  W_i=(w_1,\,\cdots,\,w_i),\quad
                  P_i=(p_1,\,\cdots,\,p_i)
   \end{array}
\end{equation}

\bigskip
\noindent Some  properties of the above recursion, see for
instance,  Abaffy and Spedicato (1989), \cite{AbSp 89}, are
listed below that are the basic formulae for  use later on.
\begin{description}
\item[\textbf{ a. }]  Implicit factorization property
                     \begin{equation}\label{F1.8}
                     A_i^TP_i = L_i
                     \end{equation}
                    with  $ L_i $ nonsingular lower triangular.
\vspace{-8pt}
\item[\textbf{ b.}]  Null space characterizations
                   \begin{equation}\label{F1.9}
                   \begin{array}{c}
                   {\cal N} (H_{i+1}) ={\cal  R} (A_i^T),\quad
                   {\cal N} (H_{i+1}^T) ={\cal R} (W_i),\quad  \\[6pt]
                   {\cal N} (A_i) ={\cal R} (H_{i+1}^T)
                   \end{array}
                   \end{equation}
                 where ${\cal N}$= Null and ${\cal R}$=Range.
\vspace{-8pt}
\item[\textbf{ c.}]\ \ The linear variety containing all solutions to $Ax=b$
                   consists of the vectors of the form
                  \begin{equation}\label{F1.10}
                          x = x_{t+1} + H^T_{t+1}q
                 \end{equation}
                 where $ q \in I\!\! R^n $ is arbitrary.
\end{description}

This paper are organized as follows. Section 2 gives some basic
definitions , operations that will be used below. We establish ABS
algorithm of stochastic linear equations in Section 3. In section
4, we prove that $\xi_i,\alpha_i$ have some good properties,such
as $\xi_i,\alpha_i$, are both normal distributions, and give
iterative formula of their expectation and variance. In section 5,
a example is given to show the ABS algorithm of stochastic linear
equations. In section 6 $\alpha_i$ is discussed and same results
are also given.

\section{Preliminaries }
\setcounter{equation}{0}

 Consider a class of systems of stochastic
 linear equations being of the form

\begin{equation}\label{hhx5}
A\xi(w)=\eta(w)
\end{equation}

 In this paper, we give the ABS algorithm under the condition that
$\eta(w)$ is a $m$-dimensional normal distribution. And it is
shown that the iteration step length $\alpha_{i}$ is $N(V,\pi)$
and the iteration solution $\xi_{i}$ is $N_n(U,\Sigma )$.
 the $(V,\pi)$ and $(U,\Sigma) $ are determined by
four iterative formulae. Finally, it is proven that $\xi_{i+1}$
is a solution of the first $i$ equations and the step length
$\alpha_{i}$ is discussed.

 Given a random experiment with a
sample space $\Omega$, a function $\xi$ that assigns to each
element $\omega$ in $\Omega$ one and only one real number
$\xi(\omega)=x$ is called a random variable. The space of $\xi$
is the set of real number $\{x:x=\xi(\omega) \in \Omega\}$, where
$\omega \in \Omega$ means the element $\omega $ belongs to the
set $\Omega$.

If $\xi_1,\xi_2,\ldots,\xi_n$ are the $n$ random variable over
$(\Omega,\cal A,P)$, then the vector function
$(\xi_1,\xi_2,\ldots,\xi_n)$ are $\Omega$ is called a
$n$-dimensional random vector are $(\Omega,\cal A,P)$.

 If $\xi_i(i=1,2,\ldots, n)$ are $N(0,1)$ and independent
 random variables, then $$\eta=(\eta_1,\eta_2,\ldots,\eta_m)^T
=A_{mn}\xi+\mu_{m1}$$ is called $m$-dimensional random vectors,
whose probability
 density function or distribution function is simply called
 $m$-dimensional normal distribution, denote $\eta \sim N_m(\mu,AA^T)$. It is easy to
 know that $\xi=(\xi_1,\xi_2,\ldots,\xi_n)^T  \sim N_n(0,I_n)$.\\

  \begin{definition}\label{hhx6}
It is said that the random variable $\xi$ is equal to the random
variable $\eta$, denoted by $\xi=\eta .$ if distributions of $\xi$
and $\eta$ are the same and $E\xi=E\eta$ \quad and \quad
$D\xi=D\eta.$
\end{definition}
\begin{definition}\label{hhx7}
 Consider the system of random linear equations
$A\xi=\eta$, where $\eta$ is given and $\xi$ is to be found.
$\xi$ is said to be a solution of this system of $\xi$ satisfies
this system in the sense of Definition1.
\end{definition}

If \\$$\xi=\eta+a,a\in\cal R$$ then
$$E\xi=E\eta+a,D\xi=D\eta.$$

\section{ABS-S Algorithm for solving System of Stochastic Linear Equations}
\setcounter{equation}{0}

\noindent
Consider a system of the stochastic linear equations
\begin{equation}\label{hhx8}
    A\xi=\eta
\end{equation}
where $\eta=(\eta_1,\eta_2,\cdots,\eta_m)^T$ is $m$-dimensional
stochastic vector$,A=(a_1,a_2,\cdots,a_m)^T \in R^{m,n}$

\noindent
 the ABS-S algorithm is defined, based on the basic ABS algorithm
 as follows.
\bigskip \noindent
 {\bf  ABS-S Algorithm: for solving systems of stochastic linear equations}
\begin{description}
\item[\textbf {(A1) }]  Initialization.\\[1pt]
                     Give  an arbitrary stochastic vector $ \xi_1 \in  I\!\!R^n $,  and an arbitrarily nonsingular matrix $ H_1 \in R^{n,n} $.\\
                     Set $ i = 1 $ and iflag=0.
\item[\textbf{(B1)}] Computer two quantities.
                     \\[1pt]
                      Compute
                       $$
                          \begin{array}{l}
                                   s_i = H_i a_i \\
                                   \tau_i=a_i^T \xi_i- \eta_i
                          \end{array}
                       $$
\item[\textbf{  (C1)}] Check the compatibility of the system of linear equations.\\[1pt]
                    If $ s_i \not= 0 $ then goto (D). \\
                    If $ s_i = 0 $ and $P( \tau_i = 0)=1 $ then set
                      $$
                                    \begin{array}{l}
                                         \xi_{i+1} = \xi_i \\
                                         H_{i+1} = H _i
                                  \end{array}
                     $$
                  and go to (F), the $ i $-th equation  is a linear combination of the previous equations. Otherwise  stop, the
                  system has no solution.
\item[\textbf{ (D1)}] Compute the search vector $p_i\in I\!\!R^n$
                      by
                  \begin{equation}\label{F1.1}
                        p_i=H_i^Tz_i
                   \end{equation}
                  where $z_{i}$, the  parameter of Broyden, is arbitrary
                  save that
                   \begin{equation}\label{F1.2}
                          z_i^T H_ia_i\neq 0
                   \end{equation}

\item[\textbf{ (E1)}] Update the random approximation $ \xi_i $ of a
solution to (\ref{F1.2}) by
                  \begin{equation}\label{F1.3}
                        \xi_{i+1} = \xi_i - \alpha _ip_i
                  \end{equation}
                where the stepsize $ \alpha _i $ is computed by
                  \begin{equation}\label{F1.4}
                        \alpha _i = \displaystyle\frac{\tau_i }{a_i^Tp_i}
                  \end{equation}
                  If $i=m$ stop; $\xi_{m+1}$ solves the system(\ref{hhx5}).
\item[\textbf{ (F1)}] Update the ( Abaffian ) matrix $ H_i. $
                   Compute
                  \begin{equation}\label{F1.5}
                        H_{i+1} = H_i - \displaystyle\frac{H_i a_i^T w_i H_i
                        }{w_i^T H_i a_i}
                  \end{equation}
                   where $ w_i \in I\!\!R^n $, the parameter of Abaffy, is arbitrary save for the condition
                  \begin{equation}\label{F1.6}
                              w_i^T H_i a_i = 1 \textrm{ or } \neq 0
                  \end{equation}
\item[\textbf{ (G1)}] Increment the index $ i $ by one and go to
                     (B1).
\end{description}

\section {The properties of $\alpha_i$ and $\xi_i$}
Let $\eta \sim N_m(v,I_m), i.e, \eta_i \sim
N(v_i,1),i=1,2,\ldots,m$ it is shown that the step length
$\alpha_i$ and the iterative solution $\xi_i$ obtained in the
ABS-S algorithm given above. Are both normal distributions, that is to say.\\
\bigskip
\noindent
 \begin{theorem}\label{hhx9}
  $$
  \tau_1= a_1^T\xi_1-\eta_1 \sim N(a_1^T \xi_1-v_1,1)
  $$
 $$
 \alpha_1=\displaystyle\frac{\tau_1}{a_1^Tp_1}
           \sim N(\displaystyle\frac{a_1^T\xi_1-v_1}{a_1^T
           p_1},\displaystyle\frac{1}{(a_1^T p_1)^2})
 $$
 \end{theorem}
\bigskip
\noindent
 \textbf{Proof.}\quad From  $\eta \sim N_m(V,I_m)$, we know that $\eta_i \sim
N(v_i,1),i=1,2,\ldots,m$. according to the linearity of the normal
distribution. again, we see that $\tau_1$ is normal distribution
and
  $\eta \sim N(v,I_m)$, $\eta_i \sim N(v_i,1),i=1,2,...,m$
 \begin{eqnarray*}
 E\tau_1&=& E(a_1^T\xi_1-\eta_1)\\
        &=&E(a_1^T \xi_1)-E(\eta_1)\\
        &=&a_1^T  \xi_1-v_1\\
 D\tau_1&=& D(a_1^T \xi_1-\eta_1)\\
        &=&D(a_1^T  \xi_1)+D(\eta_1)\\
        &=&1.
\end{eqnarray*}
\begin{eqnarray*}
 E\alpha_1&=&E\displaystyle\frac{\tau_1}{a_1^T p_1}\\
          &=&\displaystyle\frac{E\tau_1}{a_1^T  p_1}\\
          &=&\displaystyle\frac{a_1^T \xi_1-v_1}{a_1^T  p_1}.\\
D\alpha_1&=&D(\displaystyle\frac{\tau_1}{a_1^T  p_1})\\
         &=&\displaystyle\frac{D\tau_1}{a_1^T  p_1}\\
         &=&\displaystyle\frac{1}{(a_1^T  p_1)^2}.
\end{eqnarray*}

Therefore
\begin{eqnarray*}
  \alpha_1&=&\displaystyle\frac{\tau_1}{a_1^T p_1}\sim
              N(\displaystyle\frac{a_1^T  \xi_1-v_1}{a_1^T p_1},
               \displaystyle\frac{1}{(a_1^T p_1)^2})\\
 \tau_1&=& a_1^T  \xi_1-\eta_1 \sim N(a_1^T  \xi_1-v_1,1)
\end{eqnarray*}
The proof is completed. 
\begin{theorem}\label{hhx10}
  $$
  \xi_2 \sim N_n(\xi_1-\displaystyle\frac{a_1^T \xi_1-v_1}{a_1^T p_1} p_1,
  \displaystyle\frac{p_1  p_1^T}{(a_1^T p_1)^2})
  $$
\end{theorem}
\bigskip
\noindent\textbf{Proof.}\quad
 By using Theorem \ref{hhx9}, we have
 $$
 \alpha_1 \sim N(\displaystyle\frac{a_1^T \xi_1-v_1}{a_1^T  p_1},
   \displaystyle\frac{1}{(a_1^T p_1)^2})
 $$

 According to the linearity of the normal
distribution. Again, we see that $\xi_2=\xi_1-\alpha_1 p_1$ is
normal distribution, and
\begin{eqnarray*}
 E\xi_2&=&E(\xi_1-\alpha_1 p_1)=E(\xi_1)-E(\alpha_1 p_1)\\
     &=&\xi_1-E(\alpha_1) p_1=\xi_1-\displaystyle\frac{a_1^T \xi_1-v_1}{(a_1^T p_1)} p_1\\
D\xi_2&=&D(\xi_1-\alpha_1 p_1)\\
    &=&D(\xi_1)+D(\alpha_1  p_1)\\
    &=&p_1(D\alpha_1) p_1^T \\
    &=&\displaystyle\frac{p_1 p_1^T }{(a_1^T p_1)^2}\\
\end{eqnarray*}
Therefore
$$
  \xi_2 \sim N_n(\xi_1-\displaystyle
  \frac{a_1^T\xi_1-v_1}{a_1^T p_1} p_1,\displaystyle\frac{p_1 p_1^T}{(a_1^T p_1)^2})
$$
   The proof is completed.

\begin{theorem}\label{hhx11}
$ a_i p_j=0,\quad i<j$
\end{theorem}

\bigskip
\noindent
 \textbf{Proof.} see \cite{AbBs 84}, \cite{AbSp 89} 

\begin{theorem}\label{hhx12}
   $$
    \tau_i \sim N(a_i^T\xi_1-a_i^T\Sigma_{j=1}^{i-1}(E\alpha_j)p_j-v_i,1+a_i^T
      \Sigma_{j,k=1}^{i-1} cov(\alpha_j p_j,\alpha_k p_k)a_i),i \geq 2
  $$
  \end{theorem}

\bigskip
\noindent
 \textbf{Proof.}\quad The ABS-S algorithm gives rise to
 $$ \tau_i=a_i^T \xi_i-\eta_i$$
 $$ \xi_{i+1}=\xi_i-\alpha_i  p_i$$
 therefore
 \begin{eqnarray*}
  \tau_i&=&a_i^T \xi_i-\eta_i\\
        &=&a_i^T ( \xi_{i-1}-\alpha_{i-1} p_{i-1})-\eta_i\\
        &=&...\\
        &=&a_i^T  (\xi_1-\alpha_1  p_1-\alpha_2 p_2-...-\alpha_{i-1} p_{i-1})-\eta_i\\
        &=&a_i^T (\xi_1-\Sigma_{j=1}^{i-1}(\alpha_j  p_j))-\eta_i
\end{eqnarray*}
From the iterative process, we know that $\alpha_i$ is normal
distribution, $\eta_i \sim N(v_i,1)$ hence $\tau_i$ is normal
distribution, and
\begin{eqnarray*}
 E\tau_i&=&E(a_i^T (\xi_1-\Sigma_{j=1}^{i-1}(\alpha_j p_j))-\eta_i)\\
       &=&E(a_i^T(\xi_1-\Sigma_{j=1}^{i-1}(\alpha_j p_j))-E\eta_i\\
       &=&a_i^T  (\xi_1-\Sigma_{j=1}^{i-1}(E\alpha_j)p_j)-v_i\\
D\tau_i&=&D(a_i^T (\xi_1-\Sigma_{j=1}^{i-1}(\alpha_j p_j))-\eta_i)\\
      &=&D(a_i^T  (\xi_1-\Sigma_{j=1}^{i-1}(\alpha_jp_j)))+D\eta_i-2 cov(a_i^T
(\xi_1-\Sigma_{j=1}^{i-1}(\alpha_j  p_j)),\eta_i)\\
      &=&a_i^T D(\xi_1-\Sigma_{j=1}^{i-1}(\alpha_j p_j)) a_i+1+2 cov(a_i^T \Sigma_{j=1}^{i-1}(\alpha_j p_j),\eta_i)\\
      &=&a_i^T [\Sigma_{j=1}^{i-1}D(\alpha_j p_j))+\Sigma_{j \neq k}
       cov(\alpha_j p_j,\alpha_k p_k)]  a_i+1\\
      &=&a_i^T  [\Sigma_{j,k = 1}^{i-1}cov(\alpha_j  p_j,\alpha_k p_k)]
      a_i+1.
\end{eqnarray*}
The proof is completed.

\begin{theorem}\label{hhx13}
  $$
  \alpha_i \sim N(\displaystyle\frac{a_i^T\xi_1-a_i^T\Sigma_{j=1}^{i-1}(E\alpha_j)p_j-v_i}
       {a_i^T p_i},\displaystyle
       \frac{a_i^T[\Sigma_{j,k = 1}^{i-1}cov(\alpha_j  p_j,\alpha_k p_k)]
       a_i+1}{(a_i^T p_i)^2}),i \geq 2
  $$
  \end{theorem}

\noindent\textbf{Proof.}\quad
  According to
  $ \alpha_i=\displaystyle\frac{\tau_i}{a_i^T p_i}$
  and the linearity of the normal distribution, it is easy to
  known that the property holds.\\ The proof of theorem is
  completed. 

\begin{theorem}\label{hhx14}
  $$
  \xi_{i+1} \sim N_n(\xi_1-\Sigma_{j=1}^{j=i}(E\alpha_j) p_j,\Sigma_{j,k=1}^i
  cov(\alpha_j p_j,\alpha_k p_k))
  $$
 or
 $$
 \xi_{i+1} \sim N_n(E\xi_i-(E\alpha_i) p_i,D\xi_i+p_i
   D(\alpha_i) p_i^T-2 cov(\xi_i,\alpha_i p_i)), i \geq 2
 $$
 \end{theorem}

\noindent\textbf{Proof.}
  \begin{eqnarray*}
     \xi_i&=&\xi_{i-1}-\alpha_{i-1}  p_{i-1}\\
         &=&\xi_{i-2}-\alpha_{i-2}  p_{i-2}-\alpha_{i-1}  p_{i-1}\\
          &=& \ldots\\
     &=&\xi_1-\alpha_1  p_1-\alpha_2  p_2- \ldots -\alpha_{i-1} p_{i-1}\\
     &=&\xi_1-\Sigma_{j=1}^{j=i-1}\alpha_j p_j .
  \end{eqnarray*}
By $\alpha_j(j=1,2,\ldots ,i-1)$ is normal distribution, we have
$\xi_i$  is normal distribution, and
 \begin{eqnarray*}
  E\xi_i&=&E(\xi_1-\Sigma_{j=1}^{j=i-1}(\alpha_j) p_j) \\
      &=&E\xi_1-\Sigma_{j=1}^{j=i-1}E(\alpha_j)  p_j\\
      &=&\xi_1-\Sigma_{j=1}^{j=i-1}E(\alpha_j) p_j\\
  D\xi_i&=&D(\xi_1-\Sigma_{j=1}^{j=i-1}(\alpha_j) p_j) \\
      &=&D\xi_1-\Sigma_{j=1}^{j=i-1}D(\alpha_j)  p_j)\\
      &=&\Sigma_{j=1}^{j=i-1}D(\alpha_j  p_j)+\Sigma_{j\neq k}
      cov(\alpha_j p_j,\alpha_k p_k)\\
      &=&\Sigma_{j=1}^{j=i-1} D(\alpha_j p_j)+\Sigma_{j\neq k}
       cov(\alpha_j p_j ,\alpha_k  p_k))\\
      &=&\Sigma_{j,k=1}^{i-1} cov(\alpha_j p_j,\alpha_k p_k).
\end{eqnarray*}
The proof is completed.

\begin{theorem}\label{hhx15}
   $\xi_{i+1}$ is the solution of the fist $i$ equations
\end{theorem}

\bigskip
\noindent\textbf{Proof.}\quad
  Because of
  $$
  \xi_{i+1} \sim N_n(\xi_1-\Sigma_{j=1}^{j=i}E(\alpha_j) p_j,
     \Sigma_{j,k=1}^i
   cov(\alpha_j p_j,\alpha_k p_k))
  $$
 When $l<i$,
 \begin{eqnarray*}
  a_i^T\xi_{i+1}&=&a_i^T(\xi_i-\alpha_i p_i)\\
 &=&a_i^T\xi_i-a_i^T\displaystyle\frac{\tau_i}{a_i^T  p_i)}p_i\\
 &=&a_i^T\xi-\tau_i\\
 &=&\eta_i.\\
 a_l^T \xi_{i+1}&=&a_l^T (\xi_i-\alpha_i p_i)\\
                        &=&a_l^T \xi_i- a_l^T \alpha_i p_i)\\
                        &=&a_l^T \xi_i.
 \end{eqnarray*}
 Hence
$$
  E(a_l^T \xi_{i+1})=E(a_l^T \xi_i),D(a_l^T \xi_{i+1})=D(a_l^T \xi_i)
$$
$$
a_l ^T \xi_{i+1} \sim N(v_l,1),l \leq i.
$$
The proof is completed.

\begin{theorem}\label{hhx16}
  The solution of $A\xi=\eta$ is obtained in the
  finite steps by using ABS algorithm.
\end{theorem}

\noindent\textbf{Proof.}\quad
 According to the definition (\ref{hhx6})of the solution of the system of the random
 linear equations $A\xi=\eta$, it is easy to see that the theorem
 holds.\\
 The proof of Theorem is completed. 

\section{Example}
\setcounter{equation}{0}
 Stochastic linear equation

\begin{eqnarray}\label{hhx17}
 \left(\matrix{1&3&-1&0&2&0\cr 0&-2&4&1&0&0\cr 0&-4&1&0&-2&1}\right)
 \left(\matrix{\xi^(1)\cr \xi^(2)\cr \xi^(3)\cr \xi^(4)\cr \xi^(5)\cr \xi^(6)}\right)
=\left(\matrix{\eta_1\cr \eta_2\cr \eta_3}\right)
 \end{eqnarray}
 $\eta_1\sim N(6,1),\,\eta_2\sim N(12,1),\,\eta_3\sim N(2,1)$ are
 independent each other.Our aim is determining the distribution of
 $
 \xi_i=(\xi_i^{(1)},\xi_i^{(2)},\xi_i^{(3)},\xi_i^{(4)},\xi_i^{(5)},\xi_i^{(6)})
 $
 to such that $\xi_i$ is solution of equation,then we obtain the
 distribution of $\xi_i$ using ABS-S algorithm:

 Let

 $$\xi_1=(1,1,1,1,1,1)^T,\,H_1=I_6,\,z_i=a_i,\,\omega_i=a_i,\ i=1,2,3$$

 Then

 step 1 \, i=1

$$
    \tau_1=a_1^T\xi_1-\eta_1=5-\eta_1\sim N(-1,1)
$$

\bigskip
$$
    \alpha_1=\displaystyle\frac{5-\eta_1}{15}\sim N(-\displaystyle\frac{1}{15},
       \displaystyle\frac{1}{15^2})
$$

\bigskip
$$
  p_1=H_1^Tz_1,\,\omega_1=a_1
$$
\bigskip
$$
\xi_2=\xi_1-\alpha_1p_1=(\displaystyle\frac{10+\eta_1}{15},\displaystyle\frac{\eta_1}{5},
      \displaystyle\frac{20-\eta_1}{15},1,\displaystyle\frac{5+2\eta_1}{15},1)^T\sim
      N_6(U_1,\sigma_1)
$$
where
$$
U_1=(\displaystyle\frac{16}{15},\displaystyle\frac{18}{15},
     \displaystyle\frac{16}{15},1,\displaystyle\frac{17}{15},1)^T
$$
$$
\sigma_1=\displaystyle\frac{1}{15^2}\left(\matrix{1&3&-1&0&2&0\cr
             3&9&-3&0&6&0\cr -1&-3&1&0&-2&0\cr
               0&0&0&0&0&0\cr 2&6&-2&0&4&0\cr 0&0&0&0&0&0 }\right)
$$
\bigskip
$$
  H_2=H_1-\displaystyle\frac{H_1a_1^T\omega_1H_1}{\omega_1^TH_1a_1}
      =\left(\matrix{0&0&0&0&0&0\cr -3&1&0&0&0&0\cr 1&0&1&0&0&0\cr
                     0&0&0&1&0&0\cr -2&0&0&0&1&0\cr 0&0&0&0&0&1 }\right)
$$

step 2 \, i=2

$$
  \tau_2=a_2\xi_2-\eta_2=\displaystyle\frac{90-10\eta_1-15\eta_2}{15}\sim
  N(\displaystyle\frac{1}{3},\displaystyle\frac{13}{9})
$$

\bigskip
$$
 p_2=H_2^Tz_2=(10,-2,4,1,0,0)^T
$$

\bigskip
$$
    \alpha_2=\displaystyle\frac{\tau_2}{a_2^Tp_2}
            =\displaystyle\frac{90-10\eta_1-15\eta_2}{315}
            \sim N(\displaystyle\frac{1}{63},\displaystyle\frac{13}{3969})
$$
$
 \xi_3=\xi_2-\alpha_2p_2\\
$

\bigskip
$
      =(\displaystyle\frac{-740-121\eta_1+150\eta_2}{315},
                          \displaystyle\frac{190+43\eta_1-30\eta_2}{315},
                          \displaystyle\frac{40+19\eta_1+60\eta_2}{315},
                          \displaystyle\frac{220+10\eta_1+15\eta_2}{315},
                          \displaystyle\frac{5+2\eta_1}{15},1)^T
$

\bigskip
$$
 H_3=H_2-\displaystyle\frac{H_2a_2^T\omega_2H_2}{\omega_2^TH_2a_2}
      =\left(\matrix{0&0&0&0&0&0\cr 0&0&0&0&0&0\cr -5&2&1&0&0&0\cr
                     -1.5&0.5&0&1&0&0\cr -2&0&0&0&1&0\cr 0&0&0&0&0&1}\right)
$$

step 3 \,  i=3

$$
\tau_3=a_3\xi_3-\eta_3=\displaystyle\frac{-615-237\eta_1+180\eta_2-315\eta_3}{315}\sim
N(-1.61,1.89)
$$
$$
\alpha_3=\displaystyle\frac{615+237\eta_1-180\eta_2+315\eta_3}{630}
$$
$$
     p_3=H_3^Tz_3=(-1,2,1,0,-2,1)^T
$$
$
 \xi_4=\xi_3-\alpha_3p_3\\
$

\bigskip
$
      =(\displaystyle\frac{-865+479\eta_1+120\eta_2+315\eta_3}{630},
                          \displaystyle\frac{-850-388\eta_1+300\eta_2-630\eta_3}{630},
                          \displaystyle\frac{-535-199\eta_1+300\eta_2-315\eta_3}{630},\\
$

\bigskip
$
                          \displaystyle\frac{440+20\eta_1+30\eta_2}{630},
                          \displaystyle\frac{1440+558\eta_1-360\eta_2+630\eta_3}{630},
                          \displaystyle\frac{15-237\eta_1+180\eta_2-315\eta_3}{630})^T
$
$$
   \xi_4\sim N(U,\Sigma)
$$
where
$$
 U=(6.47,-1.33,1.97,1.46,2.74,0.195)^T
$$
$$
  \Sigma=\left(\matrix{0.525&-1.011&-0.331&0.037&0.958&-0.475\cr
                       -1.011&1.75&0.886&0.007&1.695&0.865\cr
                       -0.331&0.886&0.538&0.008&1.015&0.503\cr
                       0.037&0.007&0.008&0.003&-0.004&0.002\cr
                       0.958&1.695&1.015&-0.004&1.998&-0.99\cr
                       -0.495&0.865&0.503&0.002&-0.99&0.472}\right)
$$

It is clear that the following formulae hold
$$
   a_3^T\xi_4=\eta_3,\ a_2^T\xi_4=\eta_2,\, a_1^T\xi_4=\eta_1
$$
Then $\xi_4$ is a solution of stochastic linear equation
(\ref{hhx17}).


\section{Conclusion and Discussion}
\setcounter{equation}{0}

6.1  On the step length

Since
 $$
   \alpha_i \sim N(\displaystyle\frac{a_i^T\xi_1-v_i}
   {a_i^T p_i},\displaystyle\frac{1+a_i^T
   \Sigma_{j,k=1}^{i-1}cov(\alpha_jp_j,\alpha_kp_k)a_i}
   {(a_i^T p_i)^2}),\qquad i \geq 2
 $$
 Let
 $$
  E\alpha_i=\displaystyle\frac{a_i^T\xi_1-v_i}{a_i^T p_i}
 $$
 $$
 D\alpha_i=\displaystyle\frac{1+a_i^T
   \Sigma_{j,k=1}^{i-1}cov(\alpha_jp_j,\alpha_kp_k)a_i}
   {(a_i^T p_i)^2}
 $$
 Then
 $$
    \alpha_i \sim N(E\alpha_i,D\alpha_i)
 $$
 therefore
\begin{eqnarray}\label{hhx18}
   P(E\alpha_i-D\alpha_i \leq \alpha_i \leq E\alpha_i+D\alpha_i)= 0.6827
\end{eqnarray}
\begin{eqnarray}\label{hhx19}
   P(E\alpha_i-2D\alpha_i \leq \alpha_i \leq E\alpha_i+2D\alpha_i)=0.9545
\end{eqnarray}
\begin{eqnarray}\label{hhx20}
   P(E\alpha_i-3D\alpha_i \leq \alpha_i \leq
   E\alpha_i+3D\alpha_i)=0.9973
\end{eqnarray}

Therefore, we know that
$$
\alpha_{i} \in [E\alpha_i-3D\alpha_i,E\alpha_i+3D\alpha_i]
$$
 that is to say, $\alpha_i$ are determined by the initial vector $\xi_1$,
and the first component  $v_1$ of $\eta$ as well as $a_i^T p_i$.
In the practical applications, one can chooses arbitrary are
formula of (\ref{hhx18}) to (\ref{hhx20}) according to the need
and the given creditable
degree.\\

6.2  main results

\begin{description}
\item[1.]\quad system of stochastic linear equations (\ref{hhx2}) with
         $\eta\sim N_m(v, I_m)$ can be solved by ABS-S algorithm with
         the random solution $\xi$ under assumption.
\item[2.]\quad The solution $\xi_i$ generated by ABS-S algorithm is normal distribution.
\item[3.]\quad The step length $\alpha_i$ generated by ABS-S algorithm is normal distribution.
\item[4.]\quad Since the matrix $A$ is non-random, $H_i,p_i$ generated by
         ABS-S algorithm are non-random, thus the ABS-S algorithm to solve
           $A\xi=\eta$ have the same properties as that of the basic ABS
          algorithm to solve $Ax=b$
\end{description}
6.3 On problems are to be further studied
\begin{description}
  \item[1.]\quad The first problem are that when $\eta$ is an other
            distributions, whether $A\xi=\eta$ have solution and if so,
            whether it can be resolved by ABS-S algorithm.
 \item[2.]\quad The second problem is that when $a$ is random matrix, whether
             $A\xi=\eta$ have solution and if so,
             whether it can be resolved by ABS-S algorithm.
 \item[3.]\quad The third problem is that when $\eta$ is a function of random
          variable, problem are that whether $A\xi=\eta$ have solution and
           if so, whether it can be resolved by ABS-S algorithm.
 \end{description}

\end{document}